\documentclass{amsart}
\usepackage{amsfonts}
\usepackage{amssymb}
\usepackage{graphicx}
\usepackage{amsmath}

\newtheorem{theorem}{Theorem}

\newtheorem{corollary}[theorem]{Corollary}

\newtheorem{example}[theorem]{Example}

\newtheorem{lemma}[theorem]{Lemma}

\newtheorem{proposition}[theorem]{Proposition}
\newtheorem{remark}[theorem]{Remark}

\newenvironment{acknowledgment}[1][Acknowledgment]{\begin{trivlist}
\item[\hskip \labelsep {\bfseries #1}]}{\end{trivlist}}

\begin{document}

\title{Rings Over which Cyclics are Direct sums of Projective and CS or
Noetherian}
\author{C. J. Holston}
\address{Department of Mathematics, Ohio University, Athens, Ohio-45701, USA}
\email{holston@math.ohiou.edu}
\author{S. K. Jain}
\address{Department of Mathematics, Ohio University, Athens, Ohio-45701}
\email{jain@ohio.edu}
\author{A. Leroy}
\address{Department of Mathematics, University of Artois, France}
\email{andre.leroy@univ-artois.fr}
\dedicatory{Dedicated to Patrick F. Smith on his 65th birthday.}

\begin{abstract}

$R$ is called a right $WV$-ring if each simple right $R$-module is injective relative to 
proper cyclics. If $R$ is a right $WV$-ring, then $R$ is right uniform or a right $V$-ring. 
It is shown that for a right $WV$-ring $R$, $R$ is right noetherian if and only if each 
right cyclic module is a direct sum of a projective module and a CS or noetherian module. 
For a finitely generated module $M$ with projective socle over a $V$-ring $R$ such that 
every subfactor of $M$ is a direct sum of a projective module and a CS or noetherian 
module, we show $M=X\oplus T$, where $X$ is semisimple and $T$ is noetherian with zero 
socle. In the case that $M=R$, we get $R=S\oplus T$, where $S$ is a semisimple artinian 
ring, and $T$ is a direct sum of right noetherian simple rings with zero socle. In 
addition, if $R$ is a von Neumann regular ring, then it is semisimple artinian.

\end{abstract}

\maketitle

\section{Introduction and Preliminaries}

The question of studying homological properties on modules that guarantee the noetherian 
property dates back to the 1960s, when Bass and Papp showed that a ring is right noetherian 
iff direct sums of injective modules are injective. Since then, there has been continuous 
work on finding properties on classes of modules that guarantee the ring to be right 
noetherian (or some other finiteness condition). For instance, if each cyclic right module 
is: an injective module or a projective module \cite{GJS}, a direct sum of
an injective module and a projective module \cite{OS,S}, or a direct sum of a projective 
module and a module $Q$, where $Q$ is either injective or noetherian \cite{HR}, then the 
ring is right noetherian. It is also known that if every finitely generated right module 
is 
CS, then the ring is right noetherian \cite{HRY}. A celebrated theorem of Osofsky-Smith 
states that if every cyclic module is CS then R is a $qfd$-ring \cite{OS}. In this paper, 
we will consider rings over which every cyclic right module is a direct sum of a 
projective 
module and a CS or noetherian module.\\
\indent In Section 2, we first introduce a slight generalization of $V$-rings, which we 
call $WV$-rings (weakly $V$-rings). Recall that rings over which all simple modules are 
injective are known as $V$-rings \cite{MV}. A ring $R$ is called a right $WV$-ring if 
every simple right $R$-module is $R/A$-injective for any right ideal $A$ of $R$ such that 
$R/A\not\cong R$. Detailed study as to how $WV$-rings differ from $V$-rings is provided in 
Section 2. Indeed, if $R$ is a right $WV$-ring but not a right $V$-ring, then $R$ must be 
right uniform.\\
\indent In Section 3, we introduce the property (*) for an $R$-module $M$, and say that 
$M$ 
satisfies (*) if we can write $M=A\oplus B$, where $A$ is either a CS-module or a 
noetherian module, and $B$ is a projective module. Theorem 18 (a) shows that if $R$ is a 
$V$-ring and $M$ is a finitely generated $R$-module with projective socle such that each 
subfactor of $M$ satisfies (*), then $M=X\oplus T$, where $X$ is semisimple and $T$ is 
noetherian with zero socle. In particular, if $R$ is a $V$-ring such that each cyclic 
module satisfies (*), then $R=S\oplus T$, where $S$ is semisimple artinian and $T$ is a 
finite direct sum of simple noetherian rings with zero socle. Theorem 18 (b) shows that 
for 
a $WV$-ring $R$, $R$ is noetherian iff each cyclic $R$-module satisfies (*). The property 
(*) has been studied in \cite{PT} for finitely generated, as well as 2-generated, modules. 
The proofs of the main results depend upon a series of lemmas.\\
\indent We will say a module has $fud$ whenever it has finite uniform dimension. 
Throughout, we assume all rings are associative rings with identity, and all modules are 
right $R$-modules. Thus in our results, we shall omit the word 'right' when we want to say 
right noetherian, right WV-ring, etc. We shall use $\subset_e$ to denote an essential 
submodule, and $\subset_\oplus$ to denote a direct summand. For any undefined notation or 
terminology, we refer the reader to \cite{L}.

\section{$WV$-rings}

A ring $R$ is called a $WV$-ring if each simple $R$-module is $R/A$-
injective for any right ideal $A$ such that $R/A\not\cong R$ (i.e. $R/A$ is proper). Such 
rings need not be $V$-rings, as for example the ring $\mathbb{Z}_{p^2}$ for any prime $p$ 
is a $WV$-ring which is not a $V$-ring.
Let us remark that Wisbauer (\cite{W}, p.190) called a module $M$ co-semisimple if every module 
in the category $\sigma (M)$ is $M$-injective.  
Following Wisbauer's definition of co-semisimple modules, a right $WV$-ring is
 a ring for which every proper cyclic right module is co-semisimple.\\
\indent Let us first compare $V$-rings and $WV$-rings.

\begin{lemma} Let $R$ be a $WV$-ring, and $R/A$ and $R/B$ be proper cyclic modules such 
that $A\cap B=0$. Then $R$ is a $V$-ring.\end{lemma}

{\bf Proof:} Since $R$ is a WV-ring, any simple module is $R/A\times R/B$-injective. Since 
$R_R$ embeds in $R/A\times R/B$, any simple module is $R_R$ injective, i.e. $R$ is a V-
ring.\qed

\begin{theorem} Let $R$ be a WV-ring which is not a V-ring. Then $R$ must be uniform.
\end{theorem}

{\bf Proof:} Suppose $R$ is a WV-ring. If $R$ is of infinite uniform dimension, then $R$ 
contains a direct sum $A\oplus B$ where both $A$ and $B$ are infinite direct sums of 
nonzero right ideals. If $R/A\cong R$, then $R/A$ is projective, and hence there exists a 
right ideal $C$ of $R$ such that $R=C\oplus A$. But then the cyclic module $R/C$ is 
isomorphic to an infinite direct sum of nonzero modules, a contradiction. Thus $R/A$ is 
proper. Similarly $R/B$ is proper, and so $R$ is a V-ring by Lemma 1.\\
\indent Assume now that $u.dim(R)=n>1$ is finite. Then there exist closed uniform right 
ideals $U_i$ such that $\bigoplus_{i=1}^nU_i\subset_eR$. Now 
$u.dim(R/U_1)=n-1=u.dim(R/U_2)$, and so $R/U_1$ and $R/U_2$ are proper. Hence $R$ is a $V$-
ring by Lemma 1.\\
\indent So if $R$ is a $WV$-ring but not a $V$-ring, we must have $u.dim(R)=1$, i.e. $R$ is 
uniform.\qed\vskip.4cm

The proofs of the following propositions 3 and 5 are
straightforward and follow closely the classical ones given in Lam (\cite{L}, Lemma 3.75 p.99)   
or Wisbauer (\cite{W}, 23.1 p.190). 

\begin{proposition} Let $R$ be a ring such that $R/I$ is proper for any nonzero right ideal 
$I$. Then the following are equivalent:\vskip.4cm

\noindent(a) $R$ is a $WV$-ring.\\
(b) $J(R/I)=0$ for any nonzero right ideal $I$.\\
(c) Any nonzero right ideal $I\not=R$ is an intersection of maximal right ideals.\\
(d) If a simple $R$-module is contained in a cyclic module $R/I$ where $I\not=0$, then it 
is a direct summand of $R/I$.\\
In particular, the above statements are equivalent when $R$ is uniform or local.
\qed 
\end{proposition}  

\begin{corollary} If $R$ is a $WV$-ring, then $R/J(R)$ is a $V$-ring.\end{corollary}

{\bf Proof:} Let $J=J(R)$. We note that $R$ is a $V$-ring iff each right ideal ($\not=R$) 
is an intersection of maximal right ideals. If $R$ is a $WV$-ring which is not uniform, 
then $R$ is a $V$-ring (Theorem 2) and hence $J=0$. So the result is clear in this case.\\
\indent Thus we may assume $R$ is uniform. By Proposition 3, every nonzero right ideal 
($\not=R$) is an intersection of maximal right ideals. If $J=0$, then the zero ideal is 
also an intersection of maximal right ideals, and so $R(=R/J)$ is a $V$-ring. If $J\not=0$, 
then in $R/J$ all right ideals ($\not=R/J$) are intersections of maximals, and so $R/J$ is 
a $V$-ring.\qed

\begin{proposition} If $R$ is a $WV$-ring, then the following statements hold:\vskip.4cm

\noindent(a) If $I$ is a right ideal of $R$, then either $I^2=0$ or $I^2=I$.\\
(b) If $R$ is a domain, then $R$ is simple.\\
(c) If a nonzero right ideal $I$ of $R$ contains a nonzero two-sided ideal, then every 
simple $R$-module is $R/I$-injective.\\
(d) If $R$ is a von Neumann regular ring, then $R$ is a $V$-ring.\\
(e) If $R$ is a local ring and is not a $V$-ring, then $R$ has exactly three right ideals.\end{proposition}

{\bf Proof:} (a) If $R$ is a $V$-ring, it is well known that $I^2=I$ for every right ideal 
of $R$. Assume that $R$ is not a $V$-ring. Then $R$ is uniform (Theorem 2).\\
\indent Let $I\not=R$  be a right ideal and suppose $I^2\not=0$. By Proposition 3, both $I$ 
and $I^2$ are intersections of maximal right ideals. If $I^2\not=I$, there must exist a 
maximal right ideal $M$ such that $I^2\subseteq M$ but $I\not\subseteq M$. We thus have 
$R=I+M$ and we can write $1=x+m$ for some $x\in I, m\in M$. This gives $I=(x+m)I\subseteq 
xI+mI\subseteq I^2+M=M$, a contradiction. Hence $I^2=I$.\\
\indent(b) Let $0\not=a\in R$. Since $R$ is a domain, $(aR)^2\not=0$, so part (a) gives us 
$(aR)^2=aR$, i.e. $aRaR=aR$. Since $R$ is a domain this gives that $RaR=R$.\\
\indent(c) 
Let us first remark that if $T$ is a nonzero two-sided ideal of $R$, then, 
since $(\frac{R}{T})T=0$, $R/T$ is a proper cyclic module. \\
\indent Let $I$ be a nonzero right ideal of $R$ and $T$ a nonzero two-sided ideal contained 
in $I$. Since $R$ is a $WV$-ring, any simple module is $R/T$-injective, and hence any 
simple module is injective relative to $R/I\cong\dfrac{R/T}{I/T}$.\\
\indent(d) Follows from Corollary 2.4, since $J(R)=0$.\\
\indent(e) If $I\not=0$ and $I\not=R$, then $I$ is an intersection of maximal right ideals 
(Proposition 3). So $I=J(R)$. Thus $R$ has at most three right ideals.\qed

\begin{corollary} If $R$ is a $WV$-domain, then $R$ is a $V$-domain.\qed\end{corollary}

It is known that the property of $R$ being a $V$-ring is not left-right symmetric, and 
hence neither is the property of being a $WV$-ring. In fact, the property of being a $WV$-
ring that is not a $V$-ring is not left-right symmetric either, as evidenced by the 
following example, due to Faith (\cite{F1}, page 335).

\begin{example} Let $R=\left[\begin{array}[pos]{cc}a&b\\0&\sigma(a)\\\end{array}\right]$ 
where $a,b\in\mathbb{Q}(x)$ and $\sigma$ is the $\mathbb{Q}$-endomorphism of $\mathbb{Q}
(x)$ such that $\sigma(x)=x^2$.\end{example}

In this ring there are only three left ideals and hence it is a left $WV$-ring. It cannot 
be a left $V$-ring because $J(R)\not=0$. This ring is local and is not right noetherian and 
thus it cannot be a right $WV$-ring (Proposition 5 (e)).

\section{Cyclics being (CS or Noetherian) $\oplus$ Projective}

Recall that an $R$-module $M$ satisfies (*) if we can write $M=A\oplus B$, where $A$ is 
either a CS-module or a noetherian module, and $B$ is a projective module. It was shown in 
\cite{PT} that a ring $R$ is noetherian iff every 2-generated $R$-module satisfies (*). We 
remark that it is not sufficient to assume that every cyclic satisfies (*) in order for $R$ 
to be noetherian, as may be seen from the following example \cite{LL}.

\begin{example} Let $R$ be the ring of all formal power series
$$\{\sum a_ix^i|a_i\in F,i\in I\}$$
where $F$ is a field, and $I$ ranges over all well-ordered sets of nonnegative real 
numbers.\end{example}

This ring is not noetherian, but every homomorphic image is self-injective, and hence 
satisfies (*).\vskip.4cm

We begin with a result that is used throughout the paper.
We would like to thank a referee for drawing our attention to shock's result \cite{Sh}
which shortens the proof.
\begin{proposition} Let $C$ be a cyclic $R$-module such that each cyclic subfactor of $C$ 
satisfies (*), and let $S=Soc(C)$. Then $C/S$ has $fud$. Furthermore, if $R$ is a $WV$-
ring, then $C/S$ is noetherian.
\end{proposition}

{\bf Proof:} 
Let $E\subseteq_e C$ and $\frac{X}{D}$ be a cyclic subfactor of 
$\frac{C}{E}$, where $E\subseteq D \subseteq X \subseteq C$.   Then by (*),
$\frac{X}{D}=\frac{B}{D} \oplus \frac{A}{D}$ with $\frac{B}{D}$ projective and 
$\frac{A}{D}$ CS or noetherian.  Since $D$ splits from $B$, essentiality shows that 
$B=D$.   Theorem 1.3 of \cite{DHW} then applies to give that 
$\frac{C}{E}$ has $fud$.  Since $E$ was arbitrary such that $E\subseteq_e C$, this implies 
that, in particular that $\frac{C}{E}$ has $qfd$.  Then $\frac{C}{soc(C)}$ is $fud$ by 
Lemma 2.9 of \cite{DHW}.  
\noindent Now assume that $R$ is a WV-ring and let $Z\subset Y \subseteq C/E$, where, as 
above, $E\subseteq_e C$.  If $0\ne x \in rad(Y/Z)$ let $K$ be maximal in $xR$.  Since 
$\frac{C}{E}$ is singular it is proper cyclic and the simple module $\frac{xR}{K}$ is 
$\frac{C}{E}$-injective, so it splits in $\frac{Y/Z}{K}$, a contradiction since 
$\frac{xR}{K}\subseteq rad(\frac{Y/Z}{K})$.  We conclude that $rad(\frac{Y}{Z})=0$.  Then 
Theorem 3.8 of \cite{Sh} implies that $\frac{C}{E}$ is noetherian.  Now, Theorem 5.15 (1)
 in \cite{DHSW} shows that $\frac{C}{soc(C)}$ is noetherian. 
\qed\vskip.4cm

For the convenience of the reader, we state below a well-known lemma (Cf. Lemma 9.1 p. 73 
\cite{DHSW}) .

\begin{lemma} If $M$ is a finitely generated $CS$-module and $\oplus M_{i}$ is an
infinite direct sum of nonzero submodules of $M$, then $M/\oplus M_{i}$
cannot have finite uniform dimension.\qed\end{lemma}

Under a stronger assumption on a cyclic module $C$ than the condition (*), namely, if 
every cyclic
subfactor of $C$ is projective, CS, or noetherian, we show that $C$ is noetherian 
when $R$ is a $WV$-ring. This will play a key role later as we work towards the general 
result.

\begin{theorem} Let $C$ be a cyclic $R$-module such that each cyclic subfactor of $C$ is 
either CS, noetherian, or projective.\\
(a) Then $C$ has $fud$.\\
(b) If $R$ is a $WV$-ring, then $C$ is noetherian.\end{theorem}

{\bf Proof:} (a) Let $S=Soc(C)$. By Proposition 9, $C/S$ has $fud$. We show $S$ is 
finitely 
generated.\\
\indent Suppose $S$ is infinitely generated. Write $S=S_1\oplus S_2$, where $S_1$ and 
$S_2$ 
are both infinitely generated. Now by hypothesis, $C/S_1$ is either CS or noetherian or 
projective. If projective, then $S_1\subset_\oplus C$ and hence $S_1$ is cyclic, a 
contradiction as $S_1$ is infinitely generated. If noetherian, then so is $S/S_1\cong 
S_2$, 
a contradiction as $S_2$ is infinitely generated. So $C/S_1$ is CS. Furthermore, 
$(S_1\oplus S_2)/S_1\cong S_2$ is infinitely generated. Since $C/S\cong\dfrac{C/S_1}
{(S_1\oplus S_2)/S_1}$, we get a contradiction by invoking Lemma 10. Hence $S$ is finitely 
generated and so $C$ has $fud$.\\
\indent (b) Since $C/S$ is noetherian (Proposition 9), and $S$ has $fud$, it follows that 
$C$ is noetherian.\qed

\vskip.4cm

Although the proof of the next lemma is straightforward we believe that
the reader will appreciate the simple technique used in the proof.

\begin{lemma} Let $C$ be an R-module and $S=Soc(C)$. If $C/S$ is a uniform $R$-module, then for any 
two submodules $A$ and $B$ of $C$ with $A\cap B=0$, either $A$ or $B$ is 
semisimple.
\end{lemma}

{\bf Proof:} Let $K$ be a complement submodule of $A$ in $C$ containing $B$. Then $A\oplus 
K\subset_e C$. This yields $Soc(A\oplus K)=S$. Thus, $(A\oplus K)/(Soc(A\oplus K))\subseteq 
C/S$. Since $(A\oplus K)/(Soc(A\oplus K))\cong A/Soc(A)\times K/Soc(K)$ and $C/S$ is 
uniform as an $R$-module, either $A/Soc(A)$ or $K/Soc(K)$ is zero. So $A=Soc(A)$ or 
$K=Soc(K)$. In other words, either $A$ or $K$ (and hence $B$) is semisimple.\qed

\vskip.4cm

\begin{lemma} If $C$ is an $R$-module, and if $C/I=A/I\oplus B/I$ is a direct sum with 
$B/I$ a projective module, then $C=A\oplus B'$, where $B=B'\oplus I$.\end{lemma}

{\bf Proof:} From the decomposition $C/I=A/I\oplus B/I$, we have $C=A+B$, where $A\cap 
B=I$. Since $B/I$ is projective, $B=B'\oplus I$ for some $B'$. Then $C=A+(B'\oplus 
I)=A+B'$. We claim that $A\cap B'=0$.\\
\indent Let $x\in A\cap B'\subseteq A\cap B=I$. Then $x\in B'\cap I=0$. Thus $C=A\oplus 
B'$.\qed
\vskip.4cm
\begin{lemma} Let $R$ be a $WV$-ring. Let $C$ be a cyclic module with a projective socle 
(equivalently $S=Soc(C)$ is embeddable in $R$). If $C/S$ is a uniform $R$-module and each 
cyclic subfactor of $C$ satisfies (*), then $C$ is noetherian.
\end{lemma}

{\bf Proof:} First assume $R$ is a $V$-ring. Let $C'/I$ be a cyclic subfactor of $C$ and  
write $C'/I=A/I\oplus B/I$ as a direct 
sum of a CS or noetherian module and a projective module, respectively. Then, by Lemma 13, 
$C'=A\oplus B'$, where $B=B'\oplus I$. Since $\frac{C'}{Soc(C')}=\frac {C'}{C'\cap S}\cong 
\frac{C'+S}{S}\subseteq \frac{C}{S}$ is uniform, either $A$ or $B'$ is 
semisimple (Lemma 12). Note both $A$ and $B'$ are cyclic.\\
\indent Case 1: $A$ is semisimple. Since $A/I$ is semisimple cyclic and $R$ is a $V$-ring, 
$A/I$ is injective.   Moreover $A/I$ embeds in $A\subseteq Soc(C')\subseteq Soc(C)=S$.  
The hypothesis that $S$ embeds in $R$ yields that $A/I$ and hence $C'/I$ are projective.\\
\indent Case 2: $B'$ is semisimple. Since $B/I\cong B'$ is semisimple and cyclic, it is a 
finite direct sum of simple injective modules. If $A/I$ is noetherian, then clearly $C/I$ 
will be also.  
Recall that a direct sum of a CS module and a simple module is a CS module (Cf. 
Lemma 7.10 \cite{DHSW}).  Hence if $A/I$ is CS  then $A/I\oplus B/I$ is also CS.\\
\indent Thus, any cyclic subfactor of the cyclic module $C$ is either CS or noetherian or 
projective. Therefore, $C$ is noetherian by Theorem 11.\\
\indent Now if $R$ is not a $V$-ring, then it is uniform (Theorem 2). Since $S$ embeds in 
$R$, $S$ is trivially noetherian. So $C$ is noetherian because $C/S$ is noetherian 
(Proposition 9).\qed\vskip.4cm

\begin{proposition} Let $R$ be a $V$-ring. Let $C$ be a cyclic $R$-module with essential 
and projective socle. Suppose each cyclic subfactor of $C$ satisfies (*). Then $C$ is semisimple.
\end{proposition}

{\bf Proof:} Let $S=Soc(C)$. We know by Proposition 9 that $C/S$ has $fud$. Suppose 
$C/S\not=0$. Then $C/S$ contains a nonzero cyclic uniform submodule. Thus we can 
find $u\in C$ with $U=(uR+S)/S\cong uR/Soc(uR)$ uniform.  Since  
$Soc(uR)\subseteq_\oplus S$ we know that $Soc(uR)$ is projective.  
Moreover, every cyclic subfactor of $uR$ also satisfies (*).  Thus Lemma 14 implies that 
$uR$ is noetherian. 
$Soc(uR)$ is then a finite direct sum of simple modules, and hence it is injective. Since 
$Soc(uR)\subset_euR$, $uR=Soc(uR)$. This yields $U=0$, a contradiction. Thus $C/S=0$, that 
is, $C=S$, completing the proof.\qed\vskip.4cm

\begin{remark}
We note that the above proposition does not apply to a $WV$-ring $R$ which is not a 
$V$-ring. In this case $R_R$ is uniform and the only projective submodule 
of a cyclic $R$-module is the zero submodule.  
\end{remark}
\vskip.4cm

\begin{theorem} Let $R$ be a von Neumann regular $WV$-ring such that each cyclic $R$-module 
satisfies (*). Then $R$ is semisimple artinian.\end{theorem}

{\bf Proof:} By Proposition 5 (d), $R$ is a $V$-ring. Let $S=Soc(R)$. Then $R/S$ has $fud$ 
and hence it is semisimple artinian \cite{K}. Let $T$ be a complement of $S$. Then $T$ 
embeds essentially in $R/S$.  Thus $T=0$. Hence 
$S\subset_eR$. So by Proposition 15, $R$ is semisimple artinian.\qed\vskip.4cm

Finally, we prove the following general result.

\begin{theorem} (a) Let $R$ be a $V$-ring. Let $M$ be a finitely generated $R$-module with 
projective socle. Suppose each cyclic subfactor of $M$ satisfies (*). Then $M$ is 
noetherian, and 
$M=X\oplus T$, where $X$ is semisimple and $T$ is noetherian with zero socle.\\
In particular, if $R$ is a $V$-ring such that each cyclic module satisfies (*), then 
$R=S\oplus T$, where $S$ is semisimple artinian and $T$ is a finite direct sum of simple 
noetherian rings with zero socle.\\
(b) For a $WV$-ring $R$, $R$ is noetherian iff each cyclic $R$-module satisfies (*).\end{theorem}

{\bf Proof:} (a) First, assume $M$ is cyclic. Let $S_0=Soc(M)$ and let $T_0$ be a 
complement of $S_0$ in $M$. Consider the cyclic module $X_0=M/T_0$. Then $S_0$ is 
essentially embeddable in $X_0$. Since $Soc(X_0)\cong S_0$, $X_0$ is semisimple by 
Proposition 15. So $X_0$, and hence $S_0$, is a finite direct sum of simples. In particular 
$S_0$ is injective and we have 
$M=S_0\oplus T_0$. Since $M/S_0$ is noetherian (Proposition 9), $T_0$ is noetherian and 
obviously it has zero socle.\\
\indent In general, $M=\sum_{i=1}^nx_iR$. By above, each $x_iR$ is noetherian, and hence 
$M$ is noetherian. $X=Soc(M)$ is finitely generated and injective by hypothesis. Therefore 
$M=X\oplus T$, where $X$ is semisimple and $T$ is noetherian with zero socle.\\
\indent Finally, let $S=Soc(R)$ which is clearly projective in a $V$-ring and let $T$ be 
its complement. Then, as shown above, $R$ is a right noetherian $V$-ring. Therefore, $R$ 
is a direct sum of simple noetherian rings (\cite{F2}, page 70). 
So, $R=S\oplus T$, where $S$ is semisimple artinian and $T$ is a finite 
direct sum of simple noetherian rings with zero socle.\\
\indent(b) Note that if $R$ is a $WV$-ring and not a $V$-ring, then $R$ is uniform 
(Theorem 2). In this case, $Soc(R)$ is either zero or a minimal right ideal. Since 
$R/Soc(R)$ is 
noetherian (Proposition 9), we conclude that $R$ is noetherian. The converse is obvious.
\qed

\begin{remark} (a) Although Theorem 17 is a consequence of Theorem 18, the short proof 
given
for Theorem 17 is of independent interest. More generally, if $R$ is a $WV$
-ring in which each non-nil right ideal contains a nonzero idempotent and
every cyclic $R$-module satisfies (*), then $R$ is semisimple artinian \cite{K}.\\
\indent (b) Readers familiar with the Wisbauer Category $\sigma[M]$ may observe that the 
results in this paper can be more generally stated in $\sigma[M]$, where $M$ is a finitely 
generated module.\end{remark}

\begin{acknowledgment}
We would like to thank the referees for giving us valuable and helpful comments.\ 

Portions of this paper were presented at the International Conference on
Ring and Module Theory, Hacettepe University, Ankara, Turkey, August 2008,
and also at the conference on Rings and Modules, Universidade de Lisboa,
Lisbon, Portugal, September 2008.
\end{acknowledgment}


\begin{thebibliography}{12}

\bibitem{DHSW} N. V. Dung, D. V. Huynh, P.F. Smith and R. Wisbauer, 
{\it Extending modules}, Longman Scientific and technical, Harlow, 1994

\bibitem{DHW} N. V. Dung, D. V. Huynh, and R. Wisbauer, {\it On modules with finite uniform 
and Krull dimension}, Arch. Math. 57, 122-13 (1991).

\bibitem{F1} C. Faith, {\it Algebra: rings, modules and categories}, Springer-Verlag, New 
York-Berlin, 1973.

\bibitem{F2} C. Faith, {\it Rings and things and a fine array of twentieth century 
associative algebra}, Mathematical Surveys and Monographs 65, AMS, Providence, RI, 1999.

\bibitem{GJS} S. C. Goel, S. K. Jain and S. Singh, {\it Rings whose cyclic modules are 
injective or projective}, Proc. Amer. Math. Soc. 53, 16-18 (1975).

\bibitem{HR} D. V. Huynh and S. T. Rizvi, {\it An affirmative answer to a question on 
noetherian rings}, J. Algebra Appl. 7, 47-59 (2008).

\bibitem{HRY} D. V. Huynh, S. T. Rizvi and M. F. Yousif, {\it Rings whose finitely 
generated modules are extending}, J. Pure Appl. Algebra 111, 325-328 (1996).

\bibitem{K} I. Kaplansky, {\it Topological representation of algebras II}, Trans. Amer. 
Math. Soc. 68, 62-75 (1950).

\bibitem{L} T. Y. Lam, {\it Lectures on Modules and Rings}, Graduate Texts in Mathematics 
189, Springer-Verlag New York, Inc., 1999.

\bibitem{LL} L. S. Levy, {\it Commutative rings whose homomorphic images are self-
injective}, Pacific J. Math. 18, 149-153 (1966).

\bibitem{MV} G. O. Michler, O. E. Villamayor, {\it On rings whose simple modules are 
injective}, J. Algebra 25, 185-201 (1973).

\bibitem{OS} B. Osofsky and P. F. Smith, {\it Cyclic modules whose quotients have all 
complement submodules direct summands}, J. Algebra 139, 342-354 (1991).

\bibitem{PT} S. Plubtieng and H. Tansee, {\it Conditions for a ring to be noetherian or 
artinian}, Comm. Algebra 30(2), 783-786 (2002).

\bibitem{S} P. F. Smith, {\it Rings characterized by their cyclic modules}, Canad. J. Math. 
24, 93-111 (1979).

\bibitem{Sh} R.C. Shock, {\it Dual generalizations of the artinian and noetherian 
conditions}, Pacific J. Math. 54, 227-235 (1974).

\bibitem{W} R. Wisbauer, {\it Foundations of module and ring theory}, Gordon and Breach, 
Reading (19991). 

\end{thebibliography}
\end{document}